\documentclass[12pt]{article}
\usepackage[cp1251]{inputenc}
\usepackage[english,russian]{babel}
\usepackage{amsmath}
\usepackage{amssymb}
\usepackage{amsfonts}
\usepackage{graphicx}
\usepackage{xcolor}
\usepackage{comment}
\usepackage[colorlinks]{hyperref}
\usepackage{titlesec}
\usepackage{proof}

\newtheorem{lemma}{Lemma}
\newtheorem{theorem}{Theorem}

\newtheorem{proposition}{Proposition}

\begin{document}
\newcommand{\proof}{\textrm{Proof.\ }}
\newcommand{\qed}{\hfill\square}

\renewcommand{\refname}{References}
\renewcommand{\proofname}{Proof.}
\renewcommand{\figurename}{Fig.}
\thispagestyle{empty}

\begin{center}

\textbf{{\Large One can hear a discrete rectangular torus}}

\textbf{A.D. Mednykh${ }^{1,2}$, I.A. Mednykh${ }^{1,2}$, G.K. Sokolova${ }^{1,2}$}

{\small 1. Sobolev Institute of Mathematics, Novosibirsk, Russia} 

{\small 2. Novosibirsk State Unoversity, Novosibirsk, Russia}

\thanks{\sc Mednykh, A.D., Mednykh I.A., Sokolova G.K., 
One can hear a discrete rectangular torus}
\thanks{\copyright \ 2024 Mednykh, A.D., Mednykh I.A., Sokolova G.K.}
\end{center}

\begin{quote}
\noindent{\bf Abstract:} In the present paper, we prove that two discrete rectangular tori are isospectral if and only if they are isomorphic.
\medskip

\end{quote}
\bigskip

\section{Introduction}
 
In 1964, in his one-page paper ``Eigenvalues of the Laplace Operator on Some Manifolds'', J. Milnor \cite{Milnor64} described an example of two 16-dimensional flat tori that are isospectral but not isometric. Milnor's article inspired Mark Kac to write his famous paper \cite{Kac} entitled ``Can one hear the shape of a drum?''. Since these classical works, the question of what geometric properties of a manifold are determined by its Laplace operator has been the source of much research. Wolpert \cite{Wolpert} showed that a generic Riemann surface is determined by its Laplace spectrum. Nevertheless, pairs of isospectral non-isometric Riemann surfaces in every genus $\ge 4$ are known. See papers by Buser \cite{Buser86}, Brooks and Tse \cite{BrooksTse87}, and others. There are also examples of isospectral non-isometric surfaces of genus two and three with variable curvature Barden and Hyunsuk Kang \cite{BardKang}. At the same time,  R. Brooks \cite{Brooks88} proved that isospectral genus one Riemann surfaces (flat tori) are isometric. A similar result for the Klein bottle was obtained by R. Isangulov \cite{Isang2000}.

As we have seen, one-dimensional isospectral flat tori are always isometric, so a natural question is, What is the lowest dimension in which there are isospectral non-isometric flat tori? Following Milnor, the search for pairs of isospectral but non-isometric flat tori became a race towards the lowest possible dimension. Kneser found a 12-dimensional example in 1967 \cite{Kneser}. Ten years later, Kitaoka \cite{Kitaoka} reduced this to eight. In 1986, Conway and Sloane \cite{ConwaySloan} found five- and six-dimensional examples. In 1990, Schiemann \cite{Schiemann} constructed a four-dimensional example. Independently, and using a different approach, Shiota \cite{Shiota} found another example one year later in 1991. The same year, Earnest and Nipp \cite{EarnNipp} contributed one more pair.

Similar results are also known for graphs.  The problem of how graphs are defined by their spectrum is 
the subject  of many  studies (see survey by E.R.van Dam and W.H.Haemers \cite{DamHaemers}).
Peter Buser \cite{Buser92} posed an interesting problem: are two isospectral Riemann surfaces of genus two isometric? Up to our knowledge the problem is still open but, quite probably, can be solved positively. Because of the intrinsic link between Riemann surfaces and graphs, we hope that our result will be helpful to make a progress in solution of the Buser problem. More precisely, in the papers \cite{MedMedTheta} and \cite{LiuLu} the discrete counterpart for genus two Riemann surface was considered. This role is played by the so-called theta graphs. In these papers, In these papers, it was shown that two isospectral theta graphs are isomorphic.   

Recently, the following result was obtained. Consider a \textit{rectangular lattice} $\Gamma$ in $\mathbb{R}^{n}$ with diagonal basis matrix. Define a \textit{rectangular flat torus} as the associated flat torus $\mathbb{R}^{n}/\Gamma.$ Then, by (\cite{NilRowRyd}, Theorem 2.35), two rectangular flat tori are isospectral if and only if they are isometric. 

The aim of the present paper is to establish the similar result for discrete rectangular tori.

\section{Preliminary results and definitions}

\subsection{Laplacian for graphs and isospectrality} 

The graphs in this paper are undirected. Denote by $V(G)$ and $E(G),$ respectively, the set of vertices and edges of a graph $G.$ For each $u,v \in V(G),$ we set $a_{u\,v}$ to be equal to the number of edges between vertices $u$ and $v.$ The matrix $A=A(G)= \{a_{u\,v}\}_{u,v\in V(G)}$ is called the \textit{adjacency matrix} of the graph $G.$ Let $d(v)$ denote the valency of $v\in V(G)$,  and let $D=D(G)$ be the diagonal matrix indexed by $V(G)$ and with $d_{v\,v}=d(v).$

The matrix $L=L(G)=D(G)-A(G)$ is called the {\it Laplacian matrix} of $G.$ The Laplacian matrix $L(G)$ is a $n\times n$ square matrix, where $n = |V(G)|$ is the number of vertices in $G.$ Throughout the paper we shall denote by $\chi(G,x)$ the characteristic polynomial of $L(G).$ That is $\chi(G,\lambda)=\det(L(G)-\lambda I_{n}),$ where $I_{n}$ is the identity $n\times n$ matrix. For brevity, we will call $\chi(G,x)$ the \textit{Laplacian polynomial} of $G.$ The roots of $\chi(G,x)$ will be called the \textit{Laplacian eigenvalues} (or just eigenvalues) of $G.$ They are all non-negative numbers and can be arranged as follows 
$$0=\lambda_{1}(G)\le\lambda_{2}(G)\le\ldots\le\lambda_{n}(G),$$ repeated according to their multiplicity. Recall that for connected graph $G$ we always have $\lambda_{1}(G)=0$ and $\lambda_{2}(G) > 0.$

Two graphs $G$ and $H$ are called \textit{isospectral} if their Laplacian polynomials coincide $\chi(G,\lambda) = \chi(H,\lambda).$

\subsection{Cartesian product and uniqueness of decomposition}

The Cartesian product of graphs $G$ and $H$ is a graph, denoted $G\times H,$ whose vertex set is $V (G) \times V (H).$ Two vertices $(g, h)$ and $(g', h')$ are adjacent if either $g = g'$ and $h$ is adjacent to $h'$ in $H,$ or $h = h'$ and $g$ is adjacent to $g' $ in $G.$  A non-trivial graph $G$ is called {\it prime} if $G = G_1\times G_2$ implies that either $G_1$ or $G_2$ is a one vertex graph $K_1.$ It is easy to see that every non-trivial finite graph has a prime factorization with respect to the Cartesian product. The number of prime factors is at most $\log_2 |V (G)|.$

The uniqueness of the prime factor decomposition of connected graphs with respect to the Cartesian product was first shown by Sabidussi \cite{Sabidu}, and independently by Vizing \cite{Vizing}. More precisely, the following theorem holds.

\begin{theorem}[\textbf{Sabidussi-Vizing}]
A connected graph has a unique representation as a product of prime graphs, up to isomorphism and the order of the factors.
\end{theorem}

See book (\cite{ImriKlav}, Charter 6) for the proof of this theorem. 

\subsection{Theta function of graphs}

Let $G$ be a finite graph. Denote by $S_{G}$ the Laplacian spectrum of $G.$ Define $\Theta_{G}(t)$ function of graph $G$ in the following way
$$\Theta_{G}(t)= \sum\limits_{\lambda\in S_{G}}e^{-\lambda t}.$$ 

According to R. Brooks \cite{Brooks88}, we note that function $\Theta_{G}(t)$ completely defines spectrum of $G.$ Indeed, suppose that 
$$\Theta_{G}(t)=c_{1}\,e^{-\mu_{1}t}+c_{2}\,e^{-\mu_{2} t}+\ldots+c_{\ell}\,e^{-\mu_{\ell} t},\,0 \le \mu_{1}<\mu_{2} < \ldots < \mu_{\ell}.$$ Here $\mu_{j}$ are all elements of $S_{G}$ taken without repetitions and $c_{j}$ theirs respective multiplicities.

The multiplicities $c_{j}$ of eigenvalues $\mu_{j}$ are completely defined by $\Theta_{G}(t),$ because all functions $e^{-\mu_{j} t}$ are linearly independent. To find exact values of $c_{j}$ and $\mu_{j}$ we consider the expression 
$$e^{\mu t}\Theta_{G}(t)=c_{1}e^{(\mu - \mu_{1}) t} +c_{2}e^{(\mu - \mu_{2}) t} +\ldots+ 
  c_{\ell}e^{(\mu - \mu_{\ell}) t}.$$

We can conclude that $\mu_{1}=\underset{\mu}{\max}\{\underset{t\to\infty}{\lim}(e^{\mu t}\Theta_{G}(t))\}\text{ is finite}.$  Putting $\mu=\mu_{1}$ we have $\underset{t\to\infty}{\lim}e^{\mu t}\Theta_{G}(t)=c_{1}.$ 

Assume that we already found the initial values $(c_{1},\mu_{1}), (c_{2}, \mu_{2}), \ldots, (c_{k}, \mu_{k})$ for $k<\ell.$ Then $\mu_{k+1}$ is equal to $\underset{\mu}{\max}\{\underset{t\to\infty}{\lim}(e^{\mu t}(\Theta_{G}(t)-\sum\limits_{j=1}^{k}c_{j}e^{-\mu_{j} t}))\}\text{ is finite}.$ Its corresponding multiplicity $c_{k+1}$ is equal to $\underset{t\to\infty}{\lim}e^{\mu_{k+1}t}(\Theta_{G}(t)-\sum\limits_{j=1}^{k}c_{j}e^{-\mu_{j} t}).$

As a result, we get the following lemma.

\begin{lemma}\label{isolemma}
Two finite graphs $G_{1}$ and $G_{2}$ are isospectral if and only if $\Theta_{G_{1}}(t)=\Theta_{G_{2}}(t).$
\end{lemma}

Moreover, we can derive one more important property of $\Theta$-function for graphs. 

\begin{lemma}\label{thetalemma}
Let $G_1$ and $G_2$ are two finite graphs. Denote by $G_{1}\times G_{2}$ the Cartesian product of $G_{1}$ and $G_{2}.$ Then $\Theta_{G_{1}\times G_{2}}(t)=\Theta_{G_{1}}(t)\,\Theta_{G_{2}}(t).$ 
\end{lemma}

\proof Let $\lambda_{1},\lambda_{2},\ldots,\lambda_{m}$ be the spectrum of $G_{1}$ and $\mu_{1},\mu_{2},\ldots,\mu_{n}$ be the spectrum of $G_{2}$ By \cite{Mohar91} the spectrum of Cartesian product $G_{1}\times G_{2}$ is given by the set $\{\lambda_{i}+\mu_{j},\,1\le i\le m,1\le j\le n\}.$ 

Let $G=G_{1}\times G_{2}.$ We have

$$\Theta_{G}(t)= \sum\limits_{\lambda\in S_{G}}e^{-\lambda t}= \sum\limits_{i=1}^m \sum\limits_{j=1}^n e^{-(\lambda_{i}+\mu_{j}) t}=\sum\limits_{i=1}^m e^{-\lambda_{i} t}\sum\limits_{j=1}^n  e^{-\mu_{j}t}=\Theta_{G_1}(t)\,\Theta_{G_2}(t).\qed$$

\subsection{Discrete tori} We define a {\it discrete rectangular torus} as the Cartesian product $C_{m_1}\times C_{m_2}\times\ldots \times C_{m_p},$ where $C_m$ is a cyclic graph on $m$ vertices. Excluding trivial factors we assume that $m_i\ge  2,\, i=1,2,\ldots,p.$ In what follows, instead of the long term ``discrete rectangular  torus'', we will freely use the terms ``discrete torus'' or simply ``torus''.

The Cartesian product of graphs is commutative and associative in the following sense. For any two graphs 
$G_1,G_2$ the products  $G_1\times G_2$ and $G_2\times G_1$ are isomorphic. Given graphs $G_1,G_2$ and $G_3,$ the map $ ((x_1,x_2),x_3) \to (x_1,(x_2,x_3))$ is an isomorphism $(G_1 \times G_2) \times G_3 \to G_1 \times (G_2 \times G_3).$

Associativity and commutativity gives us permission to omit parentheses when dealing with products with more than two factors.This allows to present a discrete rectangular torus in the following  form
$$T=T_{m_1,\ldots,m_{p}}=C_{m_1}\times C_{m_2}\times\ldots \times C_{m_p},   \text{  where  }\,\,
1<m_1\le m_2\le\ldots\le m_p.$$

By making use of Sabidussi-Vizing  theorem, one can deduce that the latter representation of $T$ is unique.  Indeed, by Unique Square Lemma (\cite{ImriKlav}, Lemma 6.3), every nontrivial Cartesian product $G_1\times G_2$ contains a square subgraph $C_4=K_2\times K_2$. Thus, a cyclic graph $C_m$ is prime for any $m\neq4.$ 
Then, by replacing $C_4$ with two copies of $K_2$ in the above representation, we obtained an irreducible Cartesian product of prime graphs. According to the Sabidussi-Vizing theorem, it is unique. Replacing back the double product $K_2\times K_2$ with one copy of $C_4$ we have the uniqueness of the required expression.

For convenience, we will call the parameter $p$ the {\it dimension} of the  torus $T.$ As a result, we have the following important property.

\begin{lemma} Any discrete rectangular torus $T$  of dimension $p$ can be uniquely represented in the form 
$$T=C_{m_1}\times C_{m_2}\times\ldots\times C_{m_p},\,\text{ where }\,1<m_1\le m_2\le\ldots\le m_p.$$
\end{lemma}

\subsection{Spectrum of discrete torus}
We will use the following observation from \cite{Mohar91}. The Laplacian eigenvalues of the Cartesian product $G_1\times G_2$ of graphs $G_1$ and $G_2$ are equal to all the possible sums of eigenvalues of the two factors: $\lambda_{i}(G_{1})+\lambda_{j}(G_{2}),\,i = 1,\ldots , |V (G_1 )|,\,j = 1, \ldots, |V (G_2 )|.$

Let $G$ be a connected graph on $n$ vertices and 
$$0=\lambda_1 (G )<\lambda_2 (G )=a(G)\le\lambda_3 (G )\le\ldots\le\lambda_n (G )$$
be its Laplacian spectrum. Let us call the second smallest eigenvalue $a(G)$ of the spectrum an {\it algebraic connectivity} of the graph $G.$ Following to \cite{Fiedler}, from the above observation we note that $a(G_{1}\times G_{2})=\min(a(G_1),\,a(G_2)).$ 

We set $G_i=C_{m_i},$ for $i=1,\ldots,\,p.$ Then for any $i$ spectrum of $G_i$ is given by the list
$$\{ 4\sin^2(\frac{\pi j}{m_i}),\,j=0,1,\ldots, m_i-1\}.$$  
Hence, $a(G_i)=4\sin^2(\frac{\pi }{m_i}).$ Now, consider the discrete torus
$$T=T_{m_1,\ldots,m_{p}}=C_{m_1}\times C_{m_2}\times\ldots \times C_{m_p},\text{ where }
1<m_1\le m_2\le\ldots\le m_p.$$
Then, because of associativity of Cartesian product, spectrum of $T$ is given by
$$\{\sum\limits_{i=1}^p 4\sin^2(\frac{\pi j_i}{m_i}),\,i=1,\ldots,p,\,j_i=0,1,\ldots, m_i-1\}.$$ 
Also,
$$a(T)=\min(a(G_1), a(G_2),\ldots, a(G_p))=4\sin^2(\frac{\pi }{m_p}).$$

\subsection{One can hear the dimension of a discrete torus} 
In this point, we prove that one can hear the dimension of a discrete torus. In other words, if two tori are isospectral, then they have the same dimension. More precisely, we get the following statement. 

\begin{proposition}\label{dimprop} 
Consider two isospectral disctere rectangular tori $T_{m_1,\ldots,m_p},\linebreak 1<m_1\le\ldots\le m_p$ and $T_{\tilde{m}_1,\ldots,\tilde{m}_{\tilde{p}}},\,1< \tilde{m}_1\le \ldots\le \tilde{m}_{\tilde{p}}.$ Then $p=\tilde{p}.$
\end{proposition} 

\proof If $p=\tilde{p},$ we have what needs to be done. Without loss of generality, we can assume that $p<\tilde{p}.$ Since the tori  $T=T_{m_1,\ldots,m_p}$ and $\tilde{T}= T_{\tilde{m}_1,\ldots,\tilde{m}_{\tilde{p}}}$
are isospectral, they share the algebraic connectivity $a(T)=a(\tilde{T}),$ which is certainly a spectral invariant. Therefore, $4\sin^2(\frac{\pi }{m_p})=4\sin^2(\frac{\pi }{\tilde{m}_{\tilde{p}}}),$ that gives us $m_p=\tilde{m}_{\tilde{p}}.$

Recall the number of vertices $m_1\cdot\ldots\cdot m_p$ of the torus $T_{m_1,\ldots,m_p}$
is also a spectral invariant.

We represent tori $T$ and $\tilde{T}$ as the Cartesian product $T=T_{m_1,\ldots,m_{p-1}}\times C_{m_p}$
and $\tilde{T}=T_{\tilde{m}_1,\ldots,\tilde{m}_{\tilde{p}-1}}\times C_{m_p}.$ Since $T$ and $\tilde{T}$ are isospectral, by Lemma~\ref{isolemma} we have $\Theta_{T}(t)=\Theta_{\tilde{T}}(t).$ By Lemma~\ref{thetalemma}, $\Theta_{T}(t)=\Theta_{A}(t)\Theta_{C_{m_p}}(t)$ and $\Theta_{\tilde{T}}(t)=\Theta_{B}(t)\Theta_{C_{m_p}}(t),$ where $A=T_{m_1,\ldots,m_{p-1}}$ and $B=T_{\tilde{m}_1,\ldots,\tilde{m}_{\tilde{p}-1}}.$ Since $\Theta_{C_{m_p}}(t)$ is strictly positive, we have $\Theta_{A}(t)=\Theta_{B}(t).$  That is, graph $A$ and $B$ are isospectral and $a(A)=a(B).$ Hence, $m_{p-1}=\tilde{m}_{\tilde{p}-1}$ and, by the same arguments, $T_{m_1,\ldots,m_{p-2}}$ and $ T_{\tilde{m}_1,\ldots,\tilde{m}_{\tilde{p}-2}}$ are isospectral. Repeating the process, we obtain that $m_{1}=\tilde{m}_{\tilde{p}-p+1}$ and tori $T_{m_1}$ and $ T_{\tilde{m}_1,\ldots,\tilde{m}_{\tilde{p}-p+1}}$ are isospectral. The latter is impossible, since the number of vertices $m_1$ of the first torus is strictly less than $\tilde{m}_1\cdot\ldots\cdot\tilde{m}_{\tilde{p}-p+1}$ of the second one.$\qed$

\section{Main results}

The main result of the paper is the following theorem.
\begin{theorem} \label{mainth}
Two discrete rectangular tori are isospectral if and only if they are isomorphic.     
\end{theorem}
\proof If two tori are isomorphic, then up to the enumeration of vertices they have the same Laplace matrix and, therefore, have the same Laplace spectrum. That is, they are isospectral.

Conversely, suppose that two discrete tori $T_{m_1,\ldots,m_p},\,1< m_1 \le\ldots\le m_p$ and 
$T_{\tilde{m}_1,\ldots,\tilde{m}_{\tilde{p}}},\,1< \tilde{m}_1\le \ldots\le \tilde{m}_{\tilde{p}}$ are isospectral. From Proposition~\ref{dimprop} we already know that $p=\tilde{p}.$ 

Now, to proof the isomorphism between $T=T_{m_1,\ldots,m_p}$ and $\tilde{T}=T_{\tilde{m}_1,\ldots,\tilde{m}_p}$ we can use induction on the dimension $p.$ For $p=1,$ isospectral tori $T_{m_1}=C_{m_1}$ and $T_{\tilde{m}_1}=C_{\tilde{m}_1}$ are cyclic graphs with the same number of vertices $m_1=\tilde{m}_1$ and, therefore, are isomorphic.  

By induction, suppose that the theorem is proved for all $k<p,$ where $p\ge 2.$ Let us prove it for $p=k.$ We  represent tori $T$ and $\tilde{T}$ as the Cartesian product  $T=A\times C_{m_p}$ and $\tilde{T}=B\times C_{\tilde{m}_p},$ where $A=T_{m_1,\ldots,m_{p-1}}$ and $B=T_{\tilde{m}_1,\ldots,\tilde{m}_{p-1}}$ are respective $(p-1)$-dimensional tori. Repeating the arguments in the proof of Proposition~\ref{dimprop}, we note that the minimal non-zero eigenvalue is common to  graphs $T$ and $\tilde{T}$ and is equal to $4\sin^2(\frac{\pi }{m_p})=4\sin^2(\frac{\pi }{\tilde{m}_p}).$ Hence, $m_p=\tilde{m}_p.$ Since $T$ and $\tilde{T}$ are isospectral, by Lemma~\ref{isolemma} we have $\Theta_{T}(t)=\Theta_{\tilde{T}}(t).$ By Lemma~\ref{thetalemma}, $\Theta_{T}(t)=\Theta_{A}(t)\Theta_{C_{m_p}}(t)$ and $\Theta_{\tilde{T}}(t)=\Theta_{B}(t)\Theta_{C_{m_p}}(t).$ Since $\Theta_{C_{m_p}}(t)>0,$ we have $\Theta_{A}(t)=\Theta_{B}(t).$  That is, graph $A$ and $B$ are isospectral. By inductive hypothesis, they are isomorphic. Then $T=A\times C_{m_p}$ and $\tilde{T}=B\times C_{m_p}$ are also isomorphic. $\qed$

\section{Conclusion}

The simplicity of the proof of Theorem~\ref{mainth} is deceptive. Such statements do not hold true even for graphs with a simple cyclic structure. For example, we can consider circulant graph $G=C_{n}(s_1,\,s_2,\cdots,s_k),$ where  $0<s_1<s_2<\ldots<s_k<n/2.$ Its Laplacian spectrum is well known. It is given by the list $S_G=\{\sum_{i=1} ^k 4\sin^2(\frac{\pi s_i j}{n}),\,j=0,1,\ldots,n-1\}.$ It was shown  by Godsil, Holton and McKay \cite{Godsil}, that graphs $C_1=C_{20}(2,3,4,7)$ and $C_2=C_{20}(3,6,7,8)$ are isospectral but not isomorphic.  
Moreover, there are no cospectral circulant graphs on less than twenty vertices.

The discrete rectangular torus considered in this paper is the simplest discrete analogue of a continuous multidimensional torus. There is a more general aspect of the discrete torus that emulates the general continuous torus coming from a multidimensional parallelepiped. See, for instance papers \cite{Louis} and \cite{LinWanZhang}. The structure of Laplacian operator as well as its spectrum is known in this case. Up to our knowledge, the problem of isospectrality of such tori has not been solved.
 
The same problem naturally arises for connection Laplacian in vector bundles on discrete tori. The properties of the connection Laplacian are the subject of investigation of many recent papers \cite{Friedli}, \cite{LinWanZhang} and others.

\end{document}